# Relationship between sectional curvature and null spaces of Lichnerowicz-type Laplacians and their smallest eigenvalues

S.E. Stepanov

All Russian institute for scientific and technical information of the Russian Academy of Science,
20, Usievicha street, 125190 Moscow, Russia
E-mail address: s.e.stepanov@mail.ru

**Abstract.** The first variant of this article contained a fatal error. Therefore, we publish second version our paper. In the present paper, we prove that the curvature operator of the second kind of a Riemannian manifold is strictly positive if its sectional curvature is strictly positive and the Ricci curvature suitably pinched. In addition, we prove several vanishing theorems for null spaces of the Lichnerowicz, Sampson, and Hodge-de Rham Laplacians and find estimates for their lowest eigenvalues on compact (without boundary) Riemannian manifolds with sectional pinched curvature.

## Introduction

Sectional curvature is one of the most important concepts in Riemannian differential geometry. The geometry and topology of manifolds with sign-definite sectional curvature have been of great interest since the early days of global Riemannian geometry and remain exciting research areas with many problems (see [1]). For example, the Hopf conjecture is well known as an open problem in global Riemannian geometry: a compact, even-dimensional Riemannian manifold with positive sectional curvature has positive Euler characteristic, while a compact $2n$-dimensional Riemannian manifold with negative sectional curvature has Euler characteristic of sign $(-1)^n$.

We are concerned in the present paper with the relationship between the sign-definite sectional curvature, the Ricci curvature suitably pinched and the curvature operator of the second kind (see [2]) of a Riemannian manifold. One of the main results of this article is following: an $n$-dimensional ($n \geq 2$) compact (without boundary) Riemannian manifold $(M, g)$ with strictly positive sectional curvature is diffeomorphic to a spherical space form if the Ricci curvature $Ric$ of $(M, g)$ satisfy the inequality $Ric < n\,\alpha$ where $\alpha$ is the minimum of the sectional curvature. Furthermore, if $(M, g)$ is a compact (without boundary) simply connected Riemannian manifold with strictly positive sectional curvature and if its Ricci curvature satisfies the non-strict inequality $Ric \leq n\,\alpha$ where $\alpha$ is the minimum of the sectional curvature, then $(M, g)$ is diffeomorphic to a Riemannian locally symmetric space. As an application of these results, we prove several vanishing theorems on null spaces of the Lichnerowicz, Sampson and Hodge-de Rham Laplacians and we find estimates of their lowest eigenvalues on $n$-dimensional complete and pinched Riemannian manifolds.

The paper is organized as follows. After § 1 devoted to definitions and statements of certain known results of the theory of Lichnerowicz-type Laplacians we consider in § 2 the null space of the ordinary Lichnerowicz Laplacian and estimate of its lowest eigenvalue. Then in §§ 3-4 we consider the null spaces of the Sampson Laplacian and Hodge-de Rham Laplacian which are examples of the Lichnerowicz-type Laplacian and estimation of their lowest eigenvalues.

## 1. A brief review of the theory of Lichnerowicz-type Laplacians

Let $(M, g)$ be a Riemannian manifold of dimension $n \geq 2$ with the Levi-Civita connection $\nabla$. We shall assume all our manifolds are smooth (i.e., $C^\infty$), complete and connected. Let $E \to M$ be a tensor bundle over $(M, g)$. The Riemannian metric $g$ and the Levi-Civita connection $\nabla$ induce a scalar product and a compatible connection on a tensor bundle $E$, respectively. In addition, the scalar products and compatible connections of $E$, as well as of $M$, will be denoted by the same symbols $g$ and $\nabla$. Moreover, we can define the $L^2$-global scalar product on $C^\infty$-sections of $E$ by the formula $\langle \theta, \vartheta \rangle = \int_M g(\theta, \vartheta) dv_g$ for $\theta, \vartheta \in C^\infty(E)$ and introduce the associated Hilbert space $L^2(E)$. Then, using the $L^2$-structures on $C^\infty(E)$, we can define the *rough Laplacian* by the formula $\bar{\Delta} = \nabla^* \circ \nabla$ where $\nabla^*$ is the formal-adjoint operator with respect to the $L^2$-global scalar product of the compatible connection $\nabla: C^\infty(E) \to C^\infty(T^*M \otimes E)$ (see [3, p. 34, 52]; [4, p. 378]).

Due to the above definitions and notations, we can consider the *Lichnerowicz-type Laplacian* $\Delta_E: C^\infty(E) \to C^\infty(E)$ defined on a tensor bundle $E$ over $(M, g)$ which satisfies the Weitzenböck decomposition formula (see [2])

$$\Delta_E = \bar{\Delta} + t\, \Re_p \qquad (1.1)$$

where $t \in \mathbb{R}$ is a suitable constant and $\Re_p$ is the *Weitzenböck curvature operator* which depends linearly on the Riemann curvature tensor $R$ and the Ricci tensor $Ric$ of $(M, g)$ with the Levi-Civita connection $\nabla$. For $t = 1$ the Laplacian $\Delta_E$ is the ordinary Lichnerowicz Laplacian (see [3, p. 54]; [5]; [6]). Various geometrical problems give rise to different values of $t \in \mathbb{R}$. For example, the theory of Lichnerowicz-type Laplacian for the case $t > 0$ can be found in the monograph [9]. In particular, the Weitzenböck curvature operator $\Re_p: \otimes^q T^*M \to \otimes^q T^*M$ satisfies the following identities (see [6, p. 315]):

$$g(\Re(\varphi), \varphi') = g(\varphi, \Re(\varphi')), \quad trace_g \Re(\varphi) = \Re(trace_g \varphi) \qquad (1.2)$$

for any $\varphi, \varphi' \in C^\infty(\otimes^q T^*M)$.

The Laplacian $\Delta_E$ is a second order elliptic linear differential operator on $C^\infty(E)$ which is symmetric with respect to the $L^2$-global scalar product. On a compact (without boundary) manifold $(M, g)$, for each nonzero $t \in \mathbb{R}$, we have the orthogonal (with respect to the $L^2$-global scalar product) decomposition formula (see [7, p. 464])

$$C^\infty(E) = Ker\, \Delta_E \oplus Im\, \Delta_E \qquad (1.3)$$

where the first component $Ker\, \Delta_E$ is the kernel of the Laplacian $\Delta_E$. Since $\Delta_E$ is elliptic, the dimension of the kernel of $\Delta_E$ is finite. In accordance with [8, p. 104], its smooth sections will be called $\Delta_E$-*harmonic*. Moreover, we define the vector space of $\Delta_E$-harmonic $C^\infty$-sections of $E \to M$ by

$$Ker \Delta_E = \{\psi \in C^\infty(E): \Delta_L \psi = 0\}$$

for $t \neq 0$. On the other hand, for $t = 0$ any smooth $\Delta_E$-*harmonic* section satisfies the equation $\nabla \varphi = 0$ on a compact (without boundary) manifold $(M, g)$. In addition, in accordance with [8, p. 104] we can define the vector space of

$$L^r(Ker\Delta_E) = \{\psi \in Ker\Delta_E : \|\psi\| \in L^r(M)\}$$

where $\|\psi\|^2 = g(\psi, \psi)$.

Next, we give three examples of different types of Lichnerowicz-type Laplacian. Everywhere in the following we denote by $\Lambda^q M$ and $S^p M$ the vector bundles of exterior differential $q$-forms ($1 \leq q \leq n-1$) and symmetric differential $p$-forms ($p \geq 2$) on $(M, g)$. Throughout this paper we will consider the vector spaces of their $C^\infty$-sections denoted by $C^\infty(\Lambda^q M)$ and $C^\infty(S^p M)$, respectively. The Riemannian metric $g$ induces a point-wise inner product on the fibres of each of these spaces.

It is well known the *divergence operator* $\delta: C^\infty(\Lambda^{q+1} M) \to C^\infty(\Lambda^q M)$ such that $\delta$ is restriction of $\nabla^*$ to the subspace $C^\infty(\Lambda^{q+1} M)$ (see [3, p. 35]). Then its formal adjoint operator $d: C^\infty(\Lambda^q M) \to C^\infty(\Lambda^{q+1} M)$ that well known as exterior derivative (see [2, p. 24]). The operator $\Delta_H = d \circ \delta + \delta \circ d: C^\infty(\Lambda^q M) \to C^\infty(\Lambda^q M)$ is the *Hodge-de Rham Laplacian* on exterior differential $q$-forms (see details in [9]). In addition, the following decomposition holds:

$$\Delta_H = \bar{\Delta} + \Re_q \qquad (1.4)$$

where $\Re_q$ is the restriction of the Weitzenböck curvature operator to exterior differential $q$-forms. Moreover, the operator $\Delta_H$ is a non-negative and self-adjoint elliptic operators with respect to the $L^2$-global scalar product on $C^\infty$-sections of $\Lambda^p M$ (see [3]; [9]).

By analogy with the above, we define the divergence operator $\delta: C^\infty(S^{p+1} M) \to C^\infty(S^p M)$ such that $\delta$ is restriction of $\nabla^*$ to the subspace $C^\infty(S^{p+1} M)$ (see [3, p. 35]). Then its formal adjoint operator $\delta^*: C^\infty(S^p M) \to C^\infty(S^{p+1} M)$ is the following: $\delta^* = (p+1) Sym \circ \nabla$ for the usual symmetry action on tensors $Sym: T^*M \otimes S^p M \to S^{p+1} M$ (see [3, p. 356]). The operator $\Delta_S = \delta \circ \delta^* + \delta^* \circ \delta: C^\infty(S^p M) \to C^\infty(S^p M)$ is the *Sampson Laplacian* defined on symmetric $p$-tensor fields (see [10] and [11]). In addition, the Weitzenböck decomposition formula for the Sampson Laplacian has the form

$$\Delta_S = \bar{\Delta} - \mathfrak{R}_p \qquad (1.5)$$

where $\mathfrak{R}_p$ is the restriction of the Weitzenböck curvature operator to symmetric $p$-tensor fields ($p \geq 2$) on $M$ (see also [10] and [11]).

At the same time, we deduce from the last equation the identities (see also [12])

$$\Delta_L = \Delta_S + 2\,\mathfrak{R}_p = \bar{\Delta} + \mathfrak{R}_p \qquad (1.6)$$

where $\Delta_L: C^\infty(S^p M) \to C^\infty(S^p M)$ is the ordinary *Lichnerowicz Laplacian*, considered in [3]; [9]; [13] and in numerous other papers and monographs. For example, the second order differential operator which defines first order deformations of Einstein metrics is the Lichnerowicz Laplacian $\Delta_S = \bar{\Delta} + \mathfrak{R}_2$ for symmetric two-tensors (see [8] and [3, Chapter 12]).

## 2. The null space of the Lichnerowicz Laplacian and estimate of its lowest eigenvalue

In this section we will consider the Lichnerowicz Laplacian $\Delta_L: C^\infty(S^p M) \to C^\infty(S^p M)$. Let $(M, g)$ be covered by a system of coordinate neighborhoods $\{U, x^1, \ldots x^n\}$, where $U$ denotes a neighborhood and $x^1, \ldots x^n$ denote local coordinates in $U$. Then we can define the natural frame $\{X_1 = \partial/\partial x^1, \ldots, X_n = \partial/\partial x^n\}$ in an arbitrary coordinate neighborhood $\{U, x^1, \ldots x^n\}$. In this case, $g_{ij} = g(X_i, X_j)$ are local components of the metric tensor $g$ with the indices $i, j, k, l, \ldots \in \{1, 2, \ldots, n\}$. We will use the same definition of the curvature tensor as in [14]. Namely, we denote by $R_{ik}$ and $R_{ikjl}$ the local components the Ricci *Ric* and curvature $R$ tensors, respectively. These components are defined by the identities (see [14, pp. 145; 203; 249])

$$\nabla_i \nabla_j \psi^k - \nabla_j \nabla_i \psi^k = \psi^m R^k{}_{mij}$$

for a $C^\infty$-vector filed $\psi = \psi^k X_k$, $R_{ijkl} = g_{im} R^m{}_{jkl}$ and $R_{kl} = R^i{}_{kil}$, where we use the Einstein summation convention. The last equations, called the Ricci identities, can be rewritten as $\nabla_i\nabla_j\psi_l - \nabla_i\nabla_j\psi_l = -\psi_l R^l{}_{mij}$ for the one-form $\psi_l = \psi^k g_{kl}$.
In the general case, the Ricci identities has the form

$$\nabla_i\nabla_j\psi_{i_1...i_p} - \nabla_i\nabla_j\psi_{i_1...i_p} = -\psi_{ki_2...i_p}R^k{}_{i_1ij} - \cdots - \psi_{i_1...i_{p-1}k}R^k{}_{i_pij}$$

for the local components $\psi_{i_1...i_p} = \psi(X_{i_1},...,X_{i_p})$ of an arbitrary $\psi \in C^\infty(\otimes^{p+1} T^*M)$. Then for the case $p = 2$, by using the Ricci identities suitably, we obtain (see also [10] and [11])

$$\Delta_S \varphi_{ij} = (\delta \circ \delta^* - \delta^* \circ \delta)\varphi_{ij} =$$
$$= -g^{kl}\nabla_k(\nabla_l\varphi_{ij} + \nabla_i\varphi_{jl} + \nabla_j\varphi_{li}) - \left(\nabla_i(-g^{kl}\nabla_k\varphi_{lj}) + \nabla_j(-g^{kl}\nabla_k\varphi_{li})\right) =$$
$$= \bar{\Delta}\,\varphi_{ij} - g^{kl}\varphi_{lj}R_{ik} - g^{kl}\varphi_{li}R_{kj} - 2R_{ikjl}\varphi^{kl} =$$
$$= \bar{\Delta}\,\varphi_{ij} - \underbrace{(R_{ik}\varphi^k_j + R_{jk}\varphi^k_i - 2R_{ikjl}\varphi^{kl})}_{\Re_2(\varphi)} \qquad (2.1)$$

for the local components $\varphi_{ij} = \varphi(X_i,X_j)$ of $\varphi \in C^\infty(S^pM)$, $\varphi^k_j = g^{ki}\varphi_{ij}$ and $\varphi^{kl} = g^{ki}g^{lj}\varphi_{ij}$ where $g^{ki}$ are the contravariant components of the metric tensor $g$. In general case, we have the formula

$$\Delta_S \varphi_{i_1...i_p} = (\delta \circ \delta^* - \delta^* \circ \delta)\varphi_{i_1...i_p} = \bar{\Delta}\,\varphi_{i_1...i_p} - \Re_p(\varphi)_{i_1...i_p} \qquad (2.2)$$

where (see [3, p. 54]; [6, p. 315])

$$\Re_p(\varphi)_{i_1...i_p} = \sum_{a=1}^{p} g^{jk}R_{i_aj}\,\varphi_{i_1...i_{a-1}ki_{a+1}...i_p} -$$
$$- \sum_{\substack{a,b=1 \\ a\neq b}}^{p} g^{jk}g^{lm}R_{i_aji_bjl}\,\varphi_{i_1...i_{a-1}ki_{a+1}...i_{b-1}li_{b-1}...i_p}$$

for the local components $\varphi_{i_1...i_p} = \varphi(X_{i_1},...,X_{i_p})$ of an arbitrary $\varphi \in C^\infty(S^pM)$.

**1.3.** Let consider the quadratic form $\mathcal{Q}_p(\varphi): S^pM \otimes S^pM \to \mathbb{R}$ defined by the equality (see also [15])

$$\mathcal{Q}_p(\varphi) = g(\Re_p(\varphi),\varphi) = \Re(\varphi)_{i_1...i_p}\varphi^{i_1...i_p} =$$
$$= p \cdot R_{ij}\varphi^{ikk_3...k_p}\varphi^j{}_{kk_3...k_p} - p(p-1)R_{ijkl}\varphi^{ikk_3...k_p}\varphi^{jl}{}_{k_3...k_p} \qquad (2.3)$$

where $\varphi^{i_1...i_q} = g^{i_1 j_1} ... g^{i_n j_n} \varphi_{j_1...j_n}$ are local contravariant components of an arbitrary $\varphi \in S^q M$. In particular, from (2.1) by direct calculation we obtain that the quadratic form $Q_2(\varphi)$ has the form

$$Q_2(\varphi) = g(\Re_2(\varphi), \varphi) = \Re(\varphi)\varphi_{ij}\varphi^{ij} =$$
$$= 2\left( R_{ij}\varphi^{ik}\varphi_k^j - R_{ijkl}\varphi^{ik}\varphi^{jl} \right). \tag{2.4}$$

In the last case, there exists an orthonormal basis $e_1, ..., e_n$ of $T_x M$ at any point $x \in M$ such that $\varphi_x(e_i, e_j) = \varepsilon_i \delta_{ij}$ for the Kronecker delta $\delta_{ij}$. Then the quadratic form $Q_2(\varphi)$ can be rewritten in the form

$$g(\Re_2(\varphi_x), \varphi_x) = \sum_{i \neq j} \sec(e_i, e_j)(\varepsilon_i - \varepsilon_j)^2 \tag{2.5}$$

where $\sec(e_i, e_j) = R(e_i, e_j, e_i, e_j)$ is the sectional curvature of $(M, g)$ in direction to $\sigma_x = span\{e_i, e_j\} \subset T_x M$ at an arbitrary point $x \in M$ (see [8, pp. 387-388]). In this case we rewrite (2.5) in the form

$$R_{ij}\varphi^{ik}\varphi_k^j - R_{ijkl}\varphi^{ik}\varphi^{jl} = \sum_{i<j} \sec(e_i, e_j)(\varepsilon_i - \varepsilon_j)^2.$$

Let $0 < \alpha \leq \sec$ and $trace_g \varphi = \varepsilon_1 + \cdots + \varepsilon_n = 0$, then it is obviously the following formulas hold:

$$\sum_{i<j} \sec(e_i, e_j)(\varepsilon_i - \varepsilon_j)^2 \geq \alpha \sum_{i<j} (\varepsilon_i - \varepsilon_j)^2 =$$

$$= \alpha \left( (n-1) \sum_i (\varepsilon_i)^2 - 2 \sum_{i<j} \varepsilon_i \varepsilon_j \right) =$$

$$= \alpha \left( n \sum_i (\varepsilon_i)^2 - \left( \sum_i \varepsilon_i \right)^2 \right) = n \alpha \|\varphi\|^2.$$

Then from the above we conclude that

$$R_{ij}\varphi^{ik}\varphi_k^j - R_{ijkl}\varphi^{ik}\varphi^{jl} \geq n \alpha \|\varphi\|^2.$$

for the local contravariant components $\varphi^{kl} = g^{ki}g^{lj}\varphi_{ij}$ of an arbitrary $\varphi \in S_0^2 M$.

**Remark.** We have considered above a smooth section of a subbundle $S_0^p M \subset S^p M$ of covariant symmetric $p$-tensor fields which are totally traceless, that is, traceless on any pair of indices, i.e., $\sum_{i=1}^n \varphi(e_i, e_i, X_3 \ldots, X_n) = 0$ for an orthonormal basis $\{e_1, e_2, \ldots, e_n\}$ for $T_x M$ at an arbitrary point $x \in M$.

The last inequality we can rewrite in the form

$$g\left(\overset{\circ}{R}(\varphi), \varphi\right) \geq n\,\alpha\,\|\varphi\|^2 - R_{ij}\varphi^{ik}\varphi_k^j. \tag{2.6}$$

We recall here that the symmetric operator $\overset{\circ}{R}: S_0^2 M \to S_0^2 M$ determined by the equations $\overset{\circ}{R}(\varphi)_{ij} = R_{iklj}\varphi^{kl}$ is called as *the curvature operator of the second kind* (see [2]). We will show below that the curvature operator of the second kind naturally arises as the term in Lichnerowicz Laplacian involving the curvature tensor (see also [11]). We say that $\overset{\circ}{R} > 0$ (resp., $\overset{\circ}{R} \geq 0$) if the eigenvalues of $\overset{\circ}{R}$ as a bilinear form on $S_0^2 M$ are positive (resp., nonnegative). In our case, we conclude from (2.6) that $\overset{\circ}{R} > 0$ (resp., $\overset{\circ}{R} \geq 0$) if $R_{ij}\varphi^{ik}\varphi_k^j < n\,\alpha\,\|\varphi\|^2$ (resp., $R_{ij}\varphi^{ik}\varphi_k^j \leq n\,\alpha\,\|\varphi\|^2$) for an arbitrary $\varphi \in S_0^2 M$. From (2.6) and the above comments, we can conclude the following statement.

**Lemma 1.** *The curvature operator of the second kind $\overset{\circ}{R}$ of an n-dimensional $(n \geq 2)$ Riemannian manifold $(M, g)$ is strictly positive (resp. nonnegative) if the sectional curvature and Ricci curvature of $(M, g)$ satisfy the compound inequality that combines two following simple inequalities: $\sec \geq \alpha$ and $\mathrm{Ric} < n\,\alpha$ (resp., $\sec \geq \alpha$ and $\mathrm{Ric} \leq n\,\alpha$) for some positive number $\alpha$.*

We recall here that by the *Hopf-Rinow-Myers theorem* (see [9, p. 251]), if a manifold $(M, g)$ is complete and satisfies $\sec \geq \alpha > 0$, then $(M, g)$ is compact and $\mathrm{diam}(M, g) \leq \pi/\sqrt{\alpha}$. At the same time, we know that if $(M, g)$ be a compact (without boundary) Riemannian manifold such that $\overset{\circ}{R}$ is strictly positive, then $M$ is diffeomorphic to a spherical space form (see [2]). At the same time, if $(M, g)$ is a compact (without boundary), simply connected Riemannian manifold with nonnegative the curvature operator of the second kind $\overset{\circ}{R}: S_0^2 M \to S_0^2 M$, then $(M, g)$ is diffeomorphic to a Riemannian locally symmetric space (see also [2]). Therefore, can formulate the following corollary (compare with other results, for example, from [18]; [19]).

**Corollary 1**. *Let $(M, g)$ be an n-dimensional $(n \geq 2)$ complete (without boundary) Riemannian manifold with the sectional curvature $\sec \geq \alpha$ for some positive number $\alpha$. If the Ricci curvature of $(M, g)$ satisfy the inequality $\mathrm{Ric} < n\,\alpha$, then $(M, g)$ is diffeomorphic to a spherical space form. Furthermore, if $(M, g)$ is a compact, simply connected Riemannian manifold with the sectional curvature $\sec \geq \alpha > 0$ and its Ricci curvature of $(M, g)$ satisfy the non-strict inequality $\mathrm{Ric} \leq n\,\alpha$, then $(M, g)$ is diffeomorphic to a Riemannian locally symmetric space.*

Corollary 1 can be rewritten in the following form.

**Corollary 2**. *An n-dimensional $(n \geq 2)$ compact (without boundary) Riemannian manifold $(M, g)$ with positive sectional curvature is diffeomorphic to a spherical space form if its Ricci curvature satisfy the inequality $\mathrm{Ric} < n\,\alpha$ where $\alpha$ is the minimum of the sectional curvature. Furthermore, if $(M, g)$ is simply connected and its Ricci curvature of $(M, g)$ satisfy the non-strict inequality $\mathrm{Ric} \leq n\,\alpha$, then $(M, g)$ is diffeomorphic to a Riemannian locally symmetric space.*

In addition, if $X \in T_x M$ is a unit vector and we can complete it to an orthonormal basis $\{X, e_2, \ldots, e_n\}$ for $T_x M$ at an arbitrary point $x \in M$, then $\mathrm{Ric}(X, X) = \sum_{i=2}^{n} \sec(X, e_i)$. Therefore, if the sectional curvature of $(M, g)$ satisfies the following inequalities: $0 < \alpha \leq \sec \leq \frac{n}{n-1}\alpha$, then the Ricci tensor $\mathrm{Ric}$ of $(M, g)$ satisfies the inequalities $0 < (n-1)\,\alpha \leq \mathrm{Ric} \leq n\,\alpha$. As a result, we conclude that $0 < (n-1)\,\alpha\,\|\varphi\|^2 \leq R_{ij}\varphi^{ik}\varphi_k^j \leq n\,\alpha\,\|\varphi\|^2$ (see [16, p. 82]). In this case, we have the inequality $\overset{\circ}{R} \geq 0$. At the same time, if $0 < \alpha \leq \sec < \frac{n}{n-1}\alpha$, then $\overset{\circ}{R} > 0$. In conclusion, we can formulate the statement.

**Corollary 3**. *The curvature operator of the second kind $\overset{\circ}{R}$ of an n-dimensional $(n \geq 2)$ Riemannian manifold $(M, g)$ is strictly positive (respectively, nonnegative)* if the sectional curvature of $(M,g)$ satisfies the double inequality $0 < \alpha \leq \sec < \frac{n}{n-1}\alpha$ (respectively, $0 < \alpha \leq \sec \leq \frac{n}{n-1}\alpha$).

In addition, we can formulate the following corollary (compare with other results, for example, from [18]; [19]).

**Corollary 4**. *An $n$-dimensional $(n \geq 2)$ complete (without boundary), simply connected Riemannian manifold $(M,g)$ is compact and diffeomorphic to a spherical space form (resp., a Riemannian locally symmetric space) if its sectional curvature satisfies the double inequality $0 < \alpha \leq \sec < \frac{n}{n-1}\alpha$ (respectively, $0 < \alpha \leq \sec \leq \frac{n}{n-1}\alpha$).*

**Remark.** We recall that a manifold $(M,g)$ has a pointwise $1/\rho$-pinched sectional curvature or, in other words, $(M,g)$ is a $\left(1 + \frac{1}{n}\right)$-pinched Riemannian manifold if $(M,g)$ has positive sectional curvature and for every point $x \in M$ the ratio of the maximum to the minimum sectional curvature at that point is less than $\rho$ (see [22]). The main result of [22] is the following: Let $(M,g)$ be a compact (without boundary) Riemannian manifold of dimension $n \geq 4$ with pointwise $1/4$-pinched sectional curvatures, then $(M,g)$ is diffeomorphic to a spherical space form. Our Corollary 3 agrees with this result.

Since the curvature operator of the second kind $\overset{\circ}{R}$ defined on traceless symmetric two-tensors on a compact (without boundary) Riemannian manifold $(M,g)$ is strictly positive, there exists a positive number $\gamma$ such that $R_{ijkl}\varphi^{il}\varphi^{jk} \geq \gamma \|\varphi\|^2$ for an arbitrary $\varphi \in S_0^2 M$ (see [20]). In this case, $\gamma$ is the smallest eigenvalue of the symmetric operator $\overset{\circ}{R}$ defined on traceless symmetric two-tensors on a compact (without boundary) manifold $(M,g)$.

At the same time, if we choose orthogonal unit vectors $X,Y \in T_x M$ at an arbitrary point $x \in M$ and define the symmetric 2-tensor field $\theta = X \otimes Y + Y \otimes X$, then $\theta \in S_0^2 M$ and by direct calculation we obtain the following equalities:

$$g\left(\overset{\circ}{R}(\theta), \theta\right) = R_{iklj}(X^i Y^j + X^j Y^i)(X^k Y^l + X^l Y^k) =$$

$$= R_{iklj}X^iY^jX^lY^k + R_{iklj}X^jY^iX^kY^l = 2\,R_{iklj}X^iY^jX^lY^k = 2\,sec(X,Y),$$

$\|\theta\|^2 = 2$ and $Ric(X) = n - 1$. Thus, from $R_{ijkl}\varphi^{il}\varphi^{jk} \geq \gamma \|\varphi\|^2$ we obtain that $sec(X,Y) \geq \gamma > 0$ for any orthogonal unit vectors $X, Y \in T_xM$ at an arbitrary point $x \in M$. At the same time, we recall that the sectional curvature of $(M, g)$ satisfy the inequality $sec \geq \alpha$ thanks to our initial assumption. Therefore, we can conclude that the double inequality $0 < \alpha \leq \gamma$ holds. In this case, the inequalities $R_{ijkl}\varphi^{il}\varphi^{jk} \geq \gamma \|\varphi\|^2 \geq \alpha \|\varphi\|^2$ are also satisfied. As a result, we have the following inequalities (see [20] and [16, pp. 82; 85])

$$R_{ijkl}\varphi^{ili_3\ldots i_p}\varphi^{jk}_{i_3\ldots i_p} \geq \gamma\,\varphi^{i_1\ldots i_p}\varphi_{i_1\ldots i_p} \geq \alpha\,\varphi^{i_1\ldots i_p}\varphi_{i_1\ldots i_p},$$

$$R_{ij}\varphi^{ik_2\ldots k_p}\varphi^j_{k_2\ldots k_p} \geq (n-1)\alpha\,\|\varphi\|^2.$$

Therefore, the quadratic form $Q_p(\varphi)$ must satisfy the following conditions:

$$Q_p(\varphi) = g(\Re_p(\varphi), \varphi) = \Re(\varphi)_{i_1\ldots i_p}\varphi^{i_1\ldots p} =$$

$$= p \cdot R_{ij}\varphi^{ikk_3\ldots k_p}\varphi^j_{kk_3\ldots k_p} - p(p-1)R_{ijkl}\varphi^{ikk_3\ldots k_p}\varphi^{jl}_{k_3\ldots k_p}$$

$$= p \cdot R_{ij}\varphi^{ikk_3\ldots k_p}\varphi^j_{kk_3\ldots k_p} + p(p-1)R_{ijkl}\varphi^{jkk_3\ldots k_p}\varphi^{il}_{k_3\ldots k_p} \geq$$

$$\geq p(n+p-2)\,\alpha\,\|\varphi\|^2 \geq 0. \qquad (2.7)$$

Equality in (2.8) is possible only if $\varphi = 0$.

Next, from (1.2) and (2.2) we conclude that $\Delta_L: C^\infty(S_0^p M) \to C^\infty(S_0^p M)$. Then by direct calculations from the Weitzenböck decomposition formula (1.6) we obtain the following formula:

$$\tfrac{1}{2}\Delta\|\varphi\|^2 = -g(\Delta_L\varphi, \varphi) + \|\nabla\varphi\|^2 + Q_p(\varphi) \qquad (2.8)$$

for $\Delta = div \circ grad$ and any $\varphi \in C^\infty(S_0^p M)$. Let $(M, g)$ be a compact (without boundary) Riemannian manifold and $\Delta_L\varphi = 0$ for some $\varphi \in C^\infty(S_0^p M)$, then from (2.7) and (2.8) we deduce the integral formula

$$0 = \int_M \left(\|\nabla\varphi\|^2 + Q_p(\varphi)\right) dv_g \geq p(n+p-2)\,\alpha\,\int_M \|\varphi\|^2\,dv_g. \qquad (2.9)$$

From the above inequality, we can conclude that there are no non-zero $\Delta_L$-harmonic traceless symmetric $p$-tensor ($p \geq 2$) fields on compact (without boundary) manifold $(M, g)$ with positive sectional curvature. Then the following theorem holds.

**Theorem 1.** *Let the sectional curvature and Ricci curvature of an $n$-dimensional $(n \geq 2)$ complete Riemannian manifold $(M, g)$ satisfy one of the following conditions:*

(i) $\quad \sec \geq \alpha$ and $Ric < n\,\alpha$,

(ii) $\quad 0 < \alpha \leq \sec < \dfrac{n}{n-1}\alpha$

*for some positive number $\alpha$, then $\dim_{\mathbb{R}} Ker\, \Delta_L = 0$ for the Lichnerowicz Laplacian $\Delta_L: C^\infty(S_0^p M) \to C^\infty(S_0^p M)$ and any $p \geq 2$.*

Next, if we suppose that $\sec \geq \alpha$ and $Ric \leq n\,\alpha$, then from the above we have the inequality $R_{ijkl}\varphi^{jk}\varphi^{il} \geq 0$ for an arbitrary $\varphi \in S_0^2 M$ with local components $\varphi_{ij}$. Therefore, the inequality $R_{ijkl}\varphi^{jkk_3\ldots k_p}\varphi^{il}_{k_3\ldots k_p} \geq 0$ is satisfied for an arbitrary $\varphi \in S_0^p M$ with local components $\varphi_{i_1\ldots i_p}$ for $p \geq 2$ (see also [7] and [20]). At the same time, we have the double inequality $(n-1)\alpha\,\|\varphi\|^2 \leq R_{ij}\varphi^{ik_2\ldots k_p}\varphi^{j}_{k_2\ldots k_p} \leq n\,\alpha\,\|\varphi\|^2$ for an arbitrary $\varphi \in S_0^p M$ with local components $\varphi_{i_1\ldots i_p}$ for $p \geq 2$ (see [16, p. 82]). As a result, from (2.3) we conclude that the quadratic form

$$Q_p(\varphi) = p \cdot R_{ij}\varphi^{ikk_3\ldots k_p}\varphi^{j}{}_{kk_3\ldots k_p} + p(p-1)R_{ijkl}\varphi^{jkk_3\ldots k_p}\varphi^{il}{}_{k_3\ldots k_p}$$

$$\geq p(n-1)\,\alpha\,\|\varphi\|^2 \geq 0. \tag{2.9}$$

Equality in (2.9) is possible only if $\varphi = 0$.

Consider the vector space $Ker\,\Delta_L$ on a compact (without boundary) Riemannian manifold $(M, g)$. If $\Delta_L \varphi = 0$ for $\varphi \in C^\infty(S_0^p M)$, then from (2.9) we deduce the integral formula

$$0 = \int_M \left(\|\nabla\varphi\|^2 + Q_p(\varphi)\right) dv_g \geq p(n-1)\,\alpha \int_M \|\varphi\|^2\, dv_g.$$

From the above inequality, we can conclude that there are no non-zero $\Delta_L$-harmonic traceless symmetric $p$-tensor $(p \geq 2)$ fields on compact (without boundary) manifold $(M, g)$ with positive sectional curvature. Then the following corollary holds.

**Corollary 3.** *Let $(M, g)$ be an $n$-dimensional ($n \geq 2$) complete Riemannian manifold and its sectional curvature and Ricci curvatures of $(M, g)$ satisfy the compound inequality that combines two following simple inequalities: $\sec \geq \alpha$ and $\mathrm{Ric} \leq n\, \alpha$ for some positive number $\alpha$, then $\dim_\mathbb{R} \mathrm{Ker}\, \Delta_L = 0$ for the Lichnerowicz Laplacian $\Delta_L : C^\infty(S_0^p M) \to C^\infty(S_0^p M)$ and for any $p \geq 2$.*

For the case of complete noncompact (without boundary) Riemannian manifolds we have

**Theorem 2.** *Let $(M, g)$ be a connected complete noncompact (without boundary) Riemannian manifold and $\Delta_L : C^\infty(S_0^p M) \to C^\infty(S_0^p M)$ for $p \geq 0$ be the Lichnerowicz Laplacian. If the sectional curvature is a non-negative at each point of $(M, g)$, then $L^r(\mathrm{Ker}\Delta_L)$ is trivial for any number $r > 1$.*

*Proof.* By direct calculation we deduce the equality $\frac{1}{2}\Delta \|\varphi\|^2 = \|\varphi\| \cdot \Delta \|\varphi\| + \|d\|\varphi\|\|^2$. Then (2.9) can be rewritten in the form

$$\|\varphi\| \cdot \Delta \|\varphi\| = -g(\Delta_L \varphi, \varphi) + Q_p(\varphi) + \|\nabla \varphi\|^2 - \|d\|\varphi\|\|^2 \quad (2.10)$$

Using the first Kato inequality $\|\nabla \varphi\|^2 \geq \|d\|\varphi\|\|^2$ (see [4, p. 380]), we obtain from (2.10) the following inequality:

$$\|\varphi\| \cdot \Delta \|\varphi\| \geq -g(\Delta_L \varphi, \varphi) + Q_p(\varphi). \quad (2.11)$$

In this case, if $\varphi$ is a $\Delta_L$-harmonic section and the quadratic form $Q_p(\varphi)$ is a non-negative at each point of $M$, then the right-hand side of the above inequality is non-negative. In this case, by the S.-T. Yau theorem the above inequality we conclude that for the any positive number $q > 1$, either $\int_M \|\varphi\|^q dv_g = \infty$ or $\|\varphi\| = C$ for some constant $C \geq 0$ (see [14]). In particular, if $\|\varphi\| \in L^q(M)$ for some the positive number $q > 1$ and the volume of $(M, g)$ is infinite, then the constant $C$ is zero. On the other hand, the fact that $\sec \geq 0$ implies negative-semidefiniteness of $Q_p(\varphi)$ for all $p \geq 2$ was previously obtained in [19]. At the same time, a complete noncompact (without boundary) Riemannian manifold with nonnegative sectional curvature has an infinite volume [21]. These two remarks complete the proof of Theorem 2.

Suppose that $(M, g)$ is an $n$-dimensional compact (without boundary) Riemannian manifold satisfying the conditions of Theorem 1. In this case, the Lichnerowicz Laplacian $\Delta_L: C^\infty(S_0^p M) \to C^\infty(S_0^p M)$ has a discrete spectrum $Spec^{(p)}\Delta_L = \{\lambda_a^{(p)}\}_{a \geq 1}$ on $(M, g)$ since it is an elliptic formally self-adjoint second order differential operator. If $\lambda_a^{(p)} \in Spec^{(p)}\Delta_L$ is a nonzero eigenvalue of $\Delta_L$ corresponding to a nonzero eigentensor $\varphi \in C^\infty(S_0^p M)$ for the case $p \geq 2$, then from (2.8) and (2.9) we obtain the integral inequality

$$\left(\lambda_a^{(p)} - p(n + p - 2)\,\alpha\right) \int_M \|\varphi\|^2 dv_g \geq 0.$$

In this case, we can conclude that the following theorem holds.

**Theorem 3.** *Let $(M, g)$ be an n-dimensional $(n \geq 2)$ complete Riemannian manifold and $\Delta_L: C^\infty(S_0^p M) \to C^\infty(S_0^p M)$ be the Lichnerowicz Laplacian. If the sectional curvature and the Ricci curvature of $(M, g)$ satisfy one of the following conditions:*

(i) $\sec \geq \alpha$ and $Ric < n\,\alpha$,

(ii) $0 < \alpha \leq \sec < \frac{n}{n-1}\alpha$

*for some positive number $\alpha$, then an arbitrary nonzero eigenvalue $\lambda_a^{(p)}$ of $\Delta_L$ satisfies the inequality $\lambda_a^{(p)} \geq p(n + p - 2)\,\alpha$.*

Finally, consider the Lichnerowitz Laplacian $\Delta_L$ acting on the vector bundle $S_0^2 M$ of symmetric traceless 2-tensor fields, which can be regarded as infinitesimal Einstein deformations of the metric $g$. Therefore, it arises in the analysis of the *stability of the Einstein metrics* (see details in [3, chapter 12]).

Namely, let $g$ be an Einstein metric on a compact (without boundary) manifold $M$, i.e., $Ric = \frac{s}{n}g$ for the scalar curvature $s$ of $(M, g)$. A symmetric 2-tensor field $\varphi$ is an *infinitesimal Einstein deformation* of $g$ if and only if it satisfies the following equation (see [3, p. 347]):

$$\Delta_L \varphi = 2\frac{s}{n} \varphi; \quad \delta\varphi = 0; \quad trace_g \varphi = 0. \tag{2.12}$$

If $\varphi \in trace_g^{-1}(0) \cap \delta^{-1}(0)$ then it is called *transverse traceless tensor* or *TT-tensor*. Therefore, if $\varphi$ is an infinitesimal Einstein deformation of $g$, then it is a *TT-tensor* and an eigenform of the Lichnerowicz Laplacian $\Delta_L$ and $2\frac{s}{n}$ is its eigenvalue. On the other hand, if $2\frac{s}{n}$ is not an eigenvalue of $\Delta_L$, then $g$ is not deformable, i. e., Einstein deformations do not exist. We also recall that if an Einstein metric $g$ does not have infinitesimal Einstein deformations, then it is called *rigid* (see [3, p. 347]). We recall here that $R_{iklj} = -R_{ikjl}$, then the Lichnerowicz Laplacian can be rewritten in the form

$$\Delta_L \varphi_{ij} = \bar{\Delta} \varphi_{ij} + R_{ik}\varphi_j^k + R_{jk}\varphi_i^k + 2 \mathring{R}(\varphi)_{ij} \tag{2.13}$$

where $\mathring{R}: S^2M \to S^2M$ is the curvature operator of the second kind.

**Remark**. In the monograph [3, p. 133], the Lichnerowicz Laplacian was rewrite in the form $\Delta_L \varphi_{ij} = \bar{\Delta} \varphi_{ij} + R_{ik}\varphi_j^k + R_{jk}\varphi_i^k - 2 \mathring{R}(\varphi)_{ij}$ where the curvature operator of the second kind was define in [3, p. 52] by the formula $\mathring{R}(\varphi)_{ij} = R_{ikjl}\varphi^{kl}$ since the local components $R_{jkli}$ of the Riemann curvature tensor $R$ were defined in the monograph by the identities $R_{jkli} = g_{im}R_{jkl}{}^m$ where $R_{jkli} = R_{kjli}$. Then using (2.12) and (2.13) we can rewrite equation (2.9) in the form

$$\frac{1}{2} \Delta \|\varphi\|^2 = \|\nabla\varphi\|^2 + 2g\left(\mathring{R}(\varphi), \varphi\right). \tag{2.14}$$

Next, from (2.14) we deduce the integral formula

$$\int_M \left(\|\nabla\varphi\|^2 + 2g\left(\mathring{R}(\varphi), \varphi\right)\right) dv_g = 0.$$

It's obvious that the inequality $\overset{\circ}{R} > 0$ conflicts with the above integral equality. In addition, the inequality $Ric < n\,\alpha$ can be rewritten in the form $s < n^2\alpha$. Therefore, the following statement holds.

**Corollary 3.** *Let $(M, g)$ be an $n$-dimensional compact (without boundary) Einstein manifold with positive sectional curvature. If its scalar curvature satisfies the inequality $s < n^2\alpha$ where $\alpha$ is the minimum of sectional curvature, then its metric is rigid.*

**Remark.** A Riemannian manifold has positive constant sectional curvature is not infinitesimal Einstein deformable (see [3, p. 132]).

## 3. The null space of the Sampson Laplacian and estimate of its lowest eigenvalue

In this section we will consider the Sampson Laplacian $\Delta_S\colon C^\infty(S^p M) \to C^\infty(S^p M)$. Let $sec \le -\beta < 0$ for some positive number $\beta$ and $trace_g \varphi = \varepsilon_1 + \cdots + \varepsilon_n = 0$, then it is obviously the following formulas hold:

$$\sum_{i<j} \sec(e_i, e_j)(\varepsilon_i - \varepsilon_j)^2 \le -\beta \sum_{i<j}(\varepsilon_i - \varepsilon_j)^2 = -n\beta \|\varphi\|^2.$$

Then from (2.4) and (2.5) we obtain the inequality

$$g\left(\overset{\circ}{R}(\varphi), \varphi\right) \le -n\beta \|\varphi\|^2 - R_{ij}\varphi^{ik}\varphi_k^j$$

for an arbitrary $\varphi \in S_0^2 M$ with local components $\varphi_{ij}$.

We say that $\overset{\circ}{R} < 0$ (resp., $\overset{\circ}{R} \le 0$) if the eigenvalues of $\overset{\circ}{R}$ as a bilinear form on $S_0^2 M$ are negative (resp., nonpositive). In our case, we conclude that $\overset{\circ}{R} < 0$ (resp., $\overset{\circ}{R} \le 0$) if $-n\beta \|\varphi\|^2 < R_{ij}\varphi^{ik}\varphi_k^j$ (resp., $-n\beta \|\varphi\|^2 \le R_{ij}\varphi^{ik}\varphi_k^j$). On the other hand, if $-\frac{n}{n-1}\beta \le sec \le -\beta < 0$, then we have $-n\beta \le Ric \le -(n-1)\beta < 0$ since $Ric(X,X) = \sum_{i=2}^n \sec(X, e_i)$ for an orthonormal basis $\{X, e_2, \ldots, e_n\}$ for $T_x M$ at any point $x \in M$. As a result, we can conclude that if $-\frac{n}{n-1}\beta < sec \le -\beta < 0$, then we have the inequality $\overset{\circ}{R} < 0$. At the same time, if $-\frac{n}{n-1}\beta \le sec \le -\beta < 0$, then we have the inequality $\overset{\circ}{R} \le 0$, respectively.

**Lemma 2.** *The curvature operator of the second kind $\overset{\circ}{R}$ of an n-dimensional $(n \geq 2)$ Riemannian manifold $(M, g)$ is negative (respect, nonpositive) if its sectional curvature satisfies the double inequality $-\frac{n}{n-1}\beta < sec \leq -\beta < 0$ (respect, $-\frac{n}{n-1}\beta \leq sec \leq -\beta < 0$).*

Since the curvature operator of the second kind $\overset{\circ}{R}$ of a compact without boundary manifold $(M, g)$ is negative, there exists a negative number $\mu$ such that $R_{ijkl}\varphi^{il}\varphi^{jk} \leq -\mu \|\varphi\|^2$ for an arbitrary $\varphi \in S_0^2 M$. In this case, $-\mu$ is the largest eigenvalue of the symmetric operator $\overset{\circ}{R}$ defined on traceless symmetric two-tensors on $(M, g)$. Furthermore, the inequalities $sec \leq -\mu < 0$ hold.

Thus, our assumption $-\frac{n}{n-1}\beta \leq sec \leq -\beta < 0$ implies that $sec \leq -\mu < 0$ for some positive number $\mu$. Therefore, we conclude that $0 < -\mu \leq -\beta$. In this case, from the last inequality and the inequality $R_{ijkl}\varphi^{il}\varphi^{jk} \leq -\mu \|\varphi\|^2$ we conclude that $R_{ijkl}\varphi^{il}\varphi^{jk} \leq -\beta \|\varphi\|^2$. In this case, the inewuality $R_{ijkl}\varphi^{jkk_3...k_p}\varphi^{il}_{k_3...k_p} \leq -\beta \|\varphi\|^2$ holds for an arbitrary $\varphi \in S_0^p M$ with local components $\varphi_{i_1...i_p}$ and for $p \geq 2$ (see also [9]; [10]). In addition, we have the following inequalities:

$$R_{ij}\varphi^{ikk_3...k_p}\varphi^{j}_{kk_3...k_p} \leq -(n-1)\mu \|\varphi\|^2 \leq -(n-1)\beta \|\varphi\|^2 < 0$$

for an arbitrary $\varphi \in S_0^p M$ with local components $\varphi_{i_1...i_p}$ for $p \geq 2$. As a result, from the above we conclude that

$$\mathcal{Q}_p(\varphi) = g(\mathfrak{R}_p(\varphi), \varphi) = \mathfrak{R}(\varphi)_{i_1...i_p}\varphi^{i_1...p} =$$
$$= p \cdot R_{ij}\varphi^{ikk_3...k_p}\varphi^{j}_{kk_3...k_p} + p(p-1)R_{ijkl}\varphi^{jkk_3...k_p}\varphi^{il}_{k_3...k_p}$$
$$\leq -p(n+p-2)\beta \|\varphi\|^2 \leq 0. \qquad (3.1)$$

Equality in (3.1) is possible only if $\varphi = 0$.

Next, from (1.2) and (2.2) we conclude that $\Delta_S: C^\infty(S_0^p M) \to C^\infty(S_0^p M)$. Then by direct calculations from the Weitzenböck decomposition formula (1.5) we obtain the following formula:

$$\frac{1}{2}\Delta \|\varphi\|^2 = -g(\Delta_S\varphi, \varphi) + \|\nabla\varphi\|^2 - \mathcal{Q}_p(\varphi) \qquad (3.4)$$

for any $\varphi \in C^\infty(S_0^p M)$. Let $(M, g)$ be a compact (without boundary) Riemannian manifold and $\Delta_S \varphi = 0$ for some $\varphi \in C^\infty(S_0^p M)$, then from (3.1) and (3.4) we deduce the inequalities

$$0 = \int_M \left( \|\nabla \varphi\|^2 - \mathcal{Q}_p(\varphi) \right) dv_g \geq p(n+p-2)\, \beta \int_M \|\varphi\|^2 dv_g \geq 0.$$

From the above, we can conclude that there are no non-zero $\Delta_S$-harmonic traceless symmetric $p$-tensor fields on compact (without boundary) $n$-dimensional ($n \geq 2$) manifold $(M, g)$ if its sectional curvature satisfies the inequalities $-\frac{n}{n-1}\beta < \sec < -\beta < 0$. Then we can formulate the following theorem.

**Theorem 4.** *If the section curvature of an $n$-dimensional ($n \geq 2$) compact (without boundary) Riemannian manifold $(M, g)$ satisfies the double inequality $-\frac{n}{n-1}\beta < \sec \leq -\beta < 0$, then $\dim_\mathbb{R} \operatorname{Ker} \Delta_L = 0$ for the Sampson Laplacian $\Delta_S \colon C^\infty(S_0^p M) \to C^\infty(S_0^p M)$.*

For the case of complete noncompact (without boundary) Riemannian manifolds we have

**Theorem 5.** *Let $(M, g)$ be a simply connected complete noncompact (without boundary) Riemannian manifold and $\Delta_S \colon C^\infty(S_0^p M) \to C^\infty(S_0^p M)$ be the Sampson Laplacian. If the sectional curvature is a non-positive at each point of $(M, g)$, then $L^r(\operatorname{Ker}\Delta_S)$ is trivial for any number $r > 1$.*

*Proof.* By direct calculation we deduce the equality

$$\frac{1}{2} \Delta \|\varphi\|^2 = \|\varphi\| \cdot \Delta \|\varphi\| + \|d\|\varphi\|\|^2.$$

Then (3.4) can be rewritten in the form

$$\|\varphi\| \cdot \Delta \|\varphi\| = - g(\Delta_L \varphi, \varphi) - \mathcal{Q}_p(\varphi) + \|\nabla \varphi\|^2 - \|d\|\varphi\|\|^2 \quad (3.5)$$

Using the first Kato inequality $\|\nabla \varphi\|^2 \geq \|d\|\varphi\|\|^2$, we obtain from (3.5) the inequality:

$$\|\varphi\| \cdot \Delta \|\varphi\| \geq - g(\Delta_S \varphi, \varphi) - \mathcal{Q}_p(\varphi). \quad (3.6)$$

In this case, if $\varphi$ is a $\Delta_S$-harmonic section and $\mathcal{Q}_p(\varphi)$ is a non-positive at each point of $M$, then the right-hand side of the above inequality is non-negative. In this case, by the S.-T. Yau theorem the above inequality we conclude that for the any positive number $r > 1$, either $\int_M \|\varphi\|^r dv_g = \infty$ or $\|\varphi\| = C$ for some constant $C \geq 0$ (see [17]; [18]). In particular, if $\|\varphi\| \in L^r(M)$ for some the positive number $r > 1$ and the volume of $(M, g)$ is infinite then the constant $C$ is zero. On the other hand, the fact that $\sec \leq 0$ implies negative-semidefiniteness of $Q_p(\varphi)$ for all $p \geq 2$ was previously obtained in [19]. At the same time, a simply connected complete noncompact (without boundary) Riemannian manifold with nonnegative sectional curvature is a *Hadamard manifold* and therefore has an infinite volume. These two remarks complete the proof of Theorem 2.

As before, suppose that $(M, g)$ is a compact (without boundary) Riemannian manifold with curvature $\sec \leq -\beta < 0$. In this case, the Sampson Laplacian $\Delta_S \colon C^\infty(S_0^p M) \to C^\infty(S_0^p M)$ has a discrete spectrum $Spec^{(p)} \Delta_L = \{\lambda_a^{(p)}\}_{a \geq 1}$ on $(M, g)$ since it is an elliptic formally self-adjoint first order differential operator. If $\lambda_a^{(p)} \in Spec^{(p)} \Delta_L$ is an eigenvalue of $\Delta_L$ corresponding to a nonzero eigentensor $\varphi \in C^\infty(S_0^p M)$ for the case $p \geq 2$, then from (2.8) and (2.9) we obtain the integral inequality

$$0 \geq \left(-\lambda_a^{(p)} + p(n + p - 2)\beta\right) \int_M \|\varphi\|^2 dv_g.$$

In this case, we can conclude that the following theorem holds.

**Theorem 3**. *Let the sectional curvature of an $n$-dimensional $(n \geq 2)$ compact (without boundary) Riemannian manifold $(M, g)$ satisfies the double inequality $-\frac{n}{n-1}\beta < \sec \leq -\beta < 0$. Then an arbitrary nonzero eigenvalue $\lambda_a^{(p)}$ of the Lichnerowicz Laplacian $\Delta_L \colon C^\infty(S_0^p M) \to C^\infty(S_0^p M)$ satisfy the inequality $\lambda_a^{(p)} \geq p(n + p - 2)\beta$.*

## 4. The null space of the Hodge-de Rham Laplacian and estimate of its lowest eigenvalue

In this section we will consider the Hodge-de Rham Laplacian $\Delta_H = C^\infty(\Lambda^q M) \to C^\infty(\Lambda^q M)$. Let $\omega_{i_1 \ldots i_p}$ are local components of an arbitrary $\omega \in C^\infty(\Lambda^q M)$. We define the symmetric traceless two-tensor $\varphi^{[i_1 i_2 \ldots i_q]}$ with components (see also [7])

$$\varphi_{jk}^{[i_1 i_2 \ldots i_q]} = \sum_{a=1}^{q} \left( \omega_{i_1 i_2 \ldots i_{a-1} j i_{a+1} i_2 \ldots i_q} g_{k i_a} + \omega_{i_1 i_2 \ldots i_{a-1} k i_{a+1} i_2 \ldots i_q} g_{j i_a} \right) - \frac{3q}{n} g_{jk} \omega_{i_1 i_2 \ldots i_q}$$

for each set of values of indices $[i_1 i_2 \ldots i_q]$ such that $1 \leq i_1 < i_2 < \cdots < i_q$. After long but simple calculations we obtain the identities (see also [7]),

$$R_{ijkl} \varphi^{il \, [i_1 i_2 \ldots i_q]} \varphi^{jk}_{[i_1 i_2 \ldots i_q]} = q \left( \frac{2(n+4q)}{n} R_{ij} \omega^{i i_2 \ldots i_q} \omega^{j}_{i_2 \ldots i_q} - \right.$$

$$\left. - 3(p-1) R_{ijkl} \omega^{ij i_3 \ldots i_q} \omega^{kl}_{i_3 \ldots i_q} - \frac{4p}{n^2} s \, \|\omega\|^2 \right);$$

$$\varphi^{ij \, [i_1 i_2 \ldots i_q]} \varphi_{ij \, [i_1 i_2 \ldots i_q]} = \frac{2q(n+2)(n-q)}{n} \|\omega\|^2.$$

Further, if we assume that $(M,g)$ is complete and one of the following two conditions holds:

(i) $\sec \geq \alpha$ and $Ric < n\alpha$,

(ii) $0 < \alpha \leq \sec < \frac{n}{n-1} \alpha$

for some positive number $\alpha$, then $(M,g)$ is compact and the inequality $\overset{\circ}{R} >$ holds. Since the curvature operator of the second kind $\overset{\circ}{R}$ of a compact (without boundary) manifold $(M,g)$ is positive, there exists a positive number $\gamma$ such that (see [20])

$$R_{ijkl} \varphi^{il \, [i_1 i_2 \ldots i_q]} \varphi^{jk}_{[i_1 i_2 \ldots i_q]} \geq \gamma \, \varphi^{ij \, [i_1 i_2 \ldots i_q]} \varphi_{ij \, [i_1 i_2 \ldots i_q]}.$$

In this case, $\gamma$ is the smallest eigenvalue of the positive symmetric curvature operator of the second kind $\overset{\circ}{R}$ defined on traceless symmetric two-tensors on $(M,g)$. At the same time, above we proved that $\alpha \leq \gamma$. Then the last inequality implies the following inequality:

$$6 \left( \frac{n+4q}{3n} R_{ij} \omega^{i i_2 \ldots i_q} \omega^{j}_{i_2 \ldots i_q} - \frac{p-1}{2} R_{ijkl} \omega^{ij i_3 \ldots i_q} \omega^{kl}_{i_3 \ldots i_q} - \frac{2q}{3n^2} s \, \|\omega\|^2 \right) \geq$$

$$\geq 6 \left( \frac{n+4q}{3n} R_{ij} \omega^{i i_2 \ldots i_q} \omega^{j}_{i_2 \ldots i_q} - \frac{p-1}{2} R_{ijkl} \omega^{ij i_3 \ldots i_q} \omega^{kl}_{i_3 \ldots i_q} - \frac{2q}{3n^2} s \, \|\omega\|^2 \right) \geq$$

$$\geq \alpha \frac{2q(n+2)(n-q)}{n} \|\omega\|^2.$$

Let consider the quadratic form $\mathcal{Q}_q(\omega)\colon \Lambda^q M \otimes \Lambda^q M \to \mathbb{R}$ defined by the equality
$$\mathcal{Q}_q(\omega) = g(\Re_q(\omega), \omega) =$$
$$= p\left( R_{ij}\omega^{ik_2\ldots k_q}\omega^{j}{}_{k_2\ldots k_q} - \frac{p-1}{2} R_{ijkl}\omega^{ikk_3\ldots k_q}\omega^{jl}{}_{k_3\ldots k_q} \right).$$

**Remark.** One can compare our definition with the definition from monograph [16, p. 53] where the authors defined similar quadratic form $F_q(\omega) = \frac{1}{q} g(\Re_q(\omega), \omega)$.

Using the above definition and notation we obtain the inequality
$$\mathcal{Q}_q(\omega) \geq q\left( \frac{2(n-2q)}{3n} R_{ij}\omega^{ii_2\ldots i_q}\omega^{j}{}_{i_2\ldots i_q} + \frac{2q}{3n^2} s \|\omega\|^2 + \alpha \frac{2q(n+2)(n-q)}{3n} \|\omega\|^2 \right),$$
where $R_{ij}\omega^{ii_2\ldots i_q}\omega^{j}{}_{i_2\ldots i_q} \geq (n-1)\gamma \|\omega\|^2 \geq (n-1)\alpha \|\omega\|^2$ and $s \geq n(n-1)\gamma \geq n(n-1)\alpha$ (see [16, pp. 82; 85]). Therefore, the last inequality equivalent to the following:
$$\mathcal{Q}_q(\omega) \geq q(n-q)\alpha \|\omega\|^2. \tag{4.1}$$

The quadratic form $\mathcal{Q}_q(\omega)$ appears in the well-known formula, valid for any exterior differential $q$-form $\omega \in C^\infty(\Lambda^q M)$:
$$\frac{1}{2} \Delta \|\varphi\|^2 = -g(\Delta_H \varphi, \varphi) + \|\nabla\varphi\|^2 + \mathcal{Q}_q(\omega). \tag{4.2}$$

Let $(M, g)$ be a compact (without boundary) Riemannian manifold and $\Delta_H \varphi = 0$ for some $\omega \in C^\infty(\Lambda^q M)$, then from (4.1) and (4.2) we deduce the integral formula
$$0 = \int_M \left( \|\nabla\varphi\|^2 + \mathcal{Q}_q(\omega) \right) dv_g \geq q(n-q)\alpha \int_M \|\varphi\|^2 dv_g.$$

Then the following theorem holds.

**Theorem 4.** *Let the sectional curvature and Ricci curvature of an n-dimensional ($n \geq 2$) complete (without boundary) Riemannian manifold $(M, g)$ satisfy one of the following conditions:*

(i)    $\sec \geq \alpha$ and $\mathrm{Ric} < n\alpha$,

(ii)    $0 < \alpha \leq \sec < \frac{n}{n-1}\alpha$

*for some positive number $\alpha$, then $\dim_{\mathbb{R}} \mathrm{Ker}\, \Delta_H = 0$ for the Hodge-de Rham Laplacian $\Delta_L\colon C^\infty(S_0^p M) \to C^\infty(S_0^p M)$ and any $p \geq 2$.*

**Remark.** Our theorem complements the classical result of Meyer (see [23]).

As before, suppose that $(M,g)$ is a complete (without boundary) Riemannian manifold and satisfies one of the following conditions holds:

(i)   $\sec \geq \alpha$ and $Ric < n\,\alpha$,

(ii)  $0 < \alpha \leq \sec < \dfrac{n}{n-1}\alpha$

for some positive number $\alpha$, then $(M,g)$ is compact and the inequality $\overset{\circ}{R} > $ holds. In this case, the Hodge-de Rham Laplacian $\Delta_H$ has a discrete spectrum $Spec^{(q)}\Delta_L = \left\{\lambda_a^{(q)}\right\}_{a\geq 1}$ on $(M,g)$ since it is an elliptic formally self-adjoint second order differential operator. If $\lambda_a^{(q)} \in Spec^{(q)}\Delta_H$ is an eigenvalue of $\Delta_H$ corresponding to a nonzero eigenform $\omega \in C^\infty(\Lambda^q M)$ for the case $q \geq 2$, then from (4.1) and (4.2) we obtain the integral inequality

$$\left(\lambda_a^{(p)} - q(n-q)\,\alpha\right)\int_M \|\varphi\|^2 dv_g \geq 0.$$

In this case, we can conclude that the following theorem holds.

**Theorem 5.** *Let $(M,g)$ be an $n$-dimensional $(n \geq 2)$ complete (without boundary) Riemannian manifold and $\Delta_H: C^\infty(\Lambda^q M) \to C^\infty(\Lambda^q M)$ be the Hodge-de Rham Laplacian. If the sectional curvature and the Ricci curvature of $(M,g)$ satisfy one of the following conditions:*

(i)   $\sec \geq \alpha$ and $Ric < n\,\alpha$,

(ii)  $0 < \alpha \leq \sec < \dfrac{n}{n-1}\alpha$

*for some positive number $\alpha$, then an arbitrary nonzero eigenvalue $\lambda_a^{(p)}$ of $\Delta_H$ satisfies the inequality $\lambda_a^{(q)} \geq q(n-q)\,\alpha$.*

**Remark.** Our theorem complements the classical result of Gallot-Meyer (see [15]).